\documentclass{amsart}
\usepackage{graphicx}
\theoremstyle{definition}

\theoremstyle{remark}

\numberwithin{equation}{section}

\begin{document}

\title{Backwards uniqueness of the mean curvature flow }%
\author{Hong Huang}%
\address{School of Mathematical Sciences, Beijing Normal University,
Laboratory of Mathematics and Complex Systems, Ministry of Education,
Beijing 100875, P.R. China}%
\email{hhuang@bnu.edu.cn}%

\subjclass{53C44}%

\keywords{mean curvature flow, backwards uniqueness, second fundamental form. }%

\begin{abstract}
In this note we prove the backwards uniqueness of the mean curvature
flow for  (codimension one) hypersurfaces in  a Euclidean  space. More precisely, let $F_t,
\widetilde{F}_t:M^n \rightarrow \mathbb{R}^{n+1}$ be two complete
solutions of the mean curvature flow on $M^n \times [0,T]$ with
bounded second fundamental forms. Suppose $F_T=\widetilde{F}_T$, then
$F_t=\widetilde{F}_t$ on $M^n \times [0,T]$. This is an analog of a
 result of Kotschwar on the Ricci flow.
\end{abstract} \maketitle

\section {Introduction}

In  [K1] Kotschwar proved backwards uniqueness of the
Ricci flow by reducing the problem to one for a suitable system of differential inequalities. Inspired by his work we prove the backwards uniqueness
of the mean curvature flow for  (codimension one) hypersurfaces in a Euclidean  space. More precisely,
we have the following

\hspace *{0.4cm}

 {\bf Theorem } Let $F_t,
\widetilde{F}_t:M^n \rightarrow \mathbb{R}^{n+1}$ be two complete
solutions of the mean curvature flow on $M^n \times [0,T]$ with
bounded second fundamental forms. Suppose $F_T=\widetilde{F}_T$, then
$F_t=\widetilde{F}_t$ on $M^n \times [0,T]$.

 \hspace *{0.4cm}

Note that the (forward) uniqueness of the mean curvature flow in any
codimension (and with more general ambient spaces) was established by Chen and Yin [CY].

As an immediate consequence of our theorem we have the following

\hspace *{0.4cm}

{\bf Corollary} Let $F_t:M^n \rightarrow \mathbb{R}^{n+1}
$ be a complete solution of the
mean curvature flow on $M^n \times [0,T]$ with bounded second
fundamental form. Let $g_t$ be the induced metric on $M^n$ via $F_t$. Suppose $\sigma$ is an isometry of $(M^n,g_T)$  such that there is a Euclidean isometry $\overline{\sigma}$ of
$\mathbb{R}^{n+1}$  satisfying $\overline{\sigma}\circ F_T=F_T
\circ \sigma$. Then there holds $\overline{\sigma}\circ F_t=F_t
\circ \sigma$ on $M^n \times [0,T]$.

\hspace *{0.4cm}

Proof of  Corollary. Note that $\overline{\sigma}\circ F_t$ and $F_t \circ
\sigma$ are two solutions to the mean curvature flow on $M^n \times [0,T]$ with bounded
second fundamental forms and with the same terminal
value, so by our theorem $\overline{\sigma}\circ F_t=F_t \circ
\sigma$ on $M^n \times [0,T]$.    \hfill{$\Box$}

\hspace *{0.4cm}

 In the next section we will give the proof of  our
theorem, which relies heavily on the methods and results in [K1] (see also [K2]).  In particular, we'll use Theorem 3.1 in [K1]. We first  reduce the proof of our theorem to that of the orientable case, so we can use the
scalar-valued second fundamental forms instead of the vector-valued forms. It is more convenient to use the scalar-valued second fundamental forms   when we do some computations  to compare  two  immersions of $M^n$ in $\mathbb{R}^{n+1}$. But towards the end of the proof we need some extra effort: We'll use the classical (Bonnet's) uniqueness theorem for hypersurfaces in a Euclidean space and Chen-Yin's uniqueness theorem for the mean curvature flow.

\section {Proof of Theorem}

To prove the Theorem we first note that we can assume that the manifold $M^n$ is connected, otherwise we can deal with each component of $M^n$.   Furthermore we can assume that $M^n$ is orientable. The reason is as follows.  If $M^n$ is not orientable, we consider the  orientation double cover $p: \hat{M}\rightarrow M$.  Let a family of immersions $F_t:M^n\rightarrow \mathbb{R}^{n+1}$ ($t\in [0,T]$) be a solution to the mean curvature
flow
\begin{displaymath}
\frac{\partial}{\partial t}F(x,t)=\vec{H}(x,t),
\end{displaymath}
where $\vec{H}(x,t)=\vec{H}_F(x,t)$ is the mean curvature vector of the immersion $F_t=F(\cdot,t)$ at the point $x\in M$.  Let $\hat{F}(\hat{x},t)=F(p(\hat{x}),t)$ for $\hat{x}\in \hat{M}$ and $t\in [0,T]$. Then
\begin{displaymath}
\frac{\partial}{\partial t}\hat{F}(\hat{x},t)=\frac{\partial}{\partial t}F(p(\hat{x}),t)=\vec{H}_F(p(\hat{x}),t)=\vec{H}_{\hat{F}}(\hat{x},t).
\end{displaymath}
That is, $\hat{F}_t=\hat{F}(\cdot,t): \hat{M} \rightarrow \mathbb{R}^{n+1}$ is also a solution to the mean curvature flow, and the proof of the Theorem in the nonorientable case  is reduced to that in the orientable case.

Now let $M^n$ be a connected, orientable, and smooth manifold, and let a family of immersions $F_t:M^n\rightarrow \mathbb{R}^{n+1}$ ($t\in [0,T]$) be a solution to the mean curvature
flow. Choose a (global) smooth, unit normal vector field  $\nu$  of the immersion $F_t$, and write  $\vec{H}=H\nu$, where $H$ is
the scalar mean curvature. Let $A=(h_{ij})$ be the (scalar) second fundamental form  of
the immersion $F_t$ w.r.t. $\nu$, $g=g_t$ be the induced metric on $M^n$ via
$F_t$, $\nabla$ be the Levi-Civita connection of $(M^n,g_t)$, and
$\Gamma^i_{jk}$ be the corresponding Christoffel symbols.  Note that $H=g^{ij}h_{ij}$, where $(g^{ij})$ is the inverse of the metric matrix $(g_{ij})$.

We have the following lemma, most of which can be
found in Huisken [H].

\hspace *{0.1cm}

{\bf Lemma 1} Along the mean curvature flow we have
\begin{eqnarray} \frac{\partial}{\partial
 t}g_{ij}&=& -2Hh_{ij}.  \\
   \frac{\partial }{\partial t}
 \Gamma^i_{jk}&=&
 -g^{il}[\nabla_j(Hh_{kl})+\nabla_k(Hh_{jl})-\nabla_l(Hh_{jk})].\\
\frac{\partial}{\partial t}h_{ij}&=& \Delta
h_{ij}-2Hh_{il}g^{lm}h_{mj}+|A|^2 h_{ij}.\\
\frac{\partial}{\partial t}\nabla_k h_{ij}&=&\Delta \nabla_k
h_{ij}+g^{pq}g^{rl}[2(h_{ki}h_{ql}-h_{kl}h_{qi})\nabla_p h_{rj}\\
  &+&2(h_{kj}h_{ql}-h_{kl}h_{qj})\nabla_p
  h_{ir}+(h_{kq}h_{pl}-h_{kl}h_{pq})\nabla_r
  h_{ij} \nonumber \\
  &+& h_{ir}\nabla_p(h_{kj}h_{ql}-h_{kl}h_{qj})+h_{rj}\nabla_p(h_{ki}h_{ql}-h_{kl}h_{qi})]\nonumber \\
  &+&g^{lm}[h_{il}(\nabla_j(Hh_{km})-\nabla_m(Hh_{kj}))+h_{lj}(\nabla_i(Hh_{km})\nonumber \\
  &-&\nabla_m(Hh_{ki}))
  -H(h_{il}\nabla_k h_{mj}+h_{jl}\nabla_k h_{mi})]\nonumber \\
  &+&\nabla_k(|A|^2
  h_{ij}).\nonumber \end{eqnarray}

\hspace *{0.1cm}

Proof. For (2.1)-(2.3) see [H]. (2.4) follows ( by a tedious
computation) from (2.2), (2.3), commutation formulas for
derivatives and the Gauss equation.   \hfill{$\Box$}

\hspace *{0.4cm}

Actually in this note we only need a rough form of the formula (2.4).

\hspace *{0.4cm}

Now let $f=g-\widetilde{g}$, $P=\nabla-\widetilde{\nabla}$,
$Q=\nabla P$,
 $S=A-\widetilde{A}$, and $U=\nabla A-\widetilde{\nabla}\widetilde{A}$, where $\widetilde{g},\widetilde{\nabla}$, etc are the
 corresponding quantities w.r.t. another family of immersions $\widetilde{F}_t:M^n\rightarrow \mathbb{R}^{n+1}$ ($t\in [0,T]$) which is also a
 solution to the mean curvature flow. Then
we have the following

\hspace *{0.4cm}

 {\bf Lemma 2}  Let $F_t$ and $\widetilde{F}_t$ be as above. We have
 \begin{eqnarray*}    \frac{\partial f}{\partial
t}&=&\widetilde{g}^{-1}\ast f  \ast \widetilde{A} \ast
\widetilde{A}+S\ast \widetilde{A}+A\ast S,\\
 \frac{\partial P}{\partial t}&=&\widetilde{g}^{-1} \ast f \ast
\widetilde{g}^{-1} \ast \widetilde{A} \ast
\widetilde{\nabla}\widetilde{A} +\widetilde{g}^{-1} \ast f \ast
\widetilde{A} \ast \widetilde{\nabla}\widetilde{A} +S \ast
\widetilde{\nabla}\widetilde{A}+A \ast  U,\\
 \frac{\partial Q}{\partial t}&=&\widetilde{g}^{-1}\ast P
\ast f \ast \widetilde{g}^{-1}\ast \widetilde{A} \ast
\widetilde{\nabla}\widetilde{A}+\widetilde{g}^{-1}\ast \widetilde{g} \ast P \ast \widetilde{g}^{-1}\ast
\widetilde{A} \ast \widetilde{\nabla}\widetilde{A}\\
&+ & \widetilde{g}^{-1}\ast f \ast \widetilde{g}^{-1} \ast
\widetilde{\nabla}\widetilde{A}\ast \widetilde{\nabla}\widetilde{A}+
\widetilde{g}^{-1}\ast f \ast \widetilde{g}^{-1} \ast \widetilde{A}
\ast {\widetilde{\nabla}}^2
\widetilde{A}\nonumber \\
&+&P\ast \widetilde{g}^{-1}\ast f \ast \widetilde{A} \ast
\widetilde{\nabla}\widetilde{A}+\widetilde{g}^{-1}\ast \widetilde{g} \ast P\ast \widetilde{A} \ast
\widetilde{\nabla}\widetilde{A}+\widetilde{g}^{-1} \ast f \ast
\widetilde{\nabla}\widetilde{A} \ast
\widetilde{\nabla}\widetilde{A}\nonumber \\
&+&\widetilde{g}^{-1} \ast f \ast \widetilde{A}   \ast
{\widetilde{\nabla}}^2 \widetilde{A}+\nabla S \ast
\widetilde{\nabla}\widetilde{A}+S \ast P \ast
\widetilde{\nabla}\widetilde{A}\nonumber \\
&+&S \ast {\widetilde{\nabla}}^2 \widetilde{A}+ \nabla A \ast U+A
\ast \nabla U + A \ast \nabla A \ast P,\nonumber  \\
(\frac{\partial}{\partial t} &-&\Delta)S = f \ast
\widetilde{g}^{-1}\ast {\widetilde{\nabla}}^2 \widetilde{A}+P \ast
\widetilde{\nabla}\widetilde{A}+Q \ast \widetilde{A}+ P \ast P \ast
\widetilde{A} \\
&+&  \widetilde{g}^{-1}\ast f \ast \widetilde{g}^{-1} \ast
\widetilde{A} \ast \widetilde{A}\ast
\widetilde{A}+\widetilde{g}^{-1}\ast f \ast  \widetilde{A} \ast
\widetilde{A}\ast \widetilde{A}+S \ast \widetilde{A}\ast
\widetilde{A}\nonumber\\
 &+&A \ast S \ast \widetilde{A}+A \ast A
\ast S,\nonumber  \\
(\frac{\partial}{\partial t}&-&\Delta)U=f \ast
\widetilde{g}^{-1}\ast {\widetilde{\nabla}}^3 \widetilde{A}+P \ast
{\widetilde{\nabla}}^2 \widetilde{A}+Q \ast
\widetilde{\nabla}\widetilde{A}+P \ast P \ast
\widetilde{\nabla}\widetilde{A}\\
&+&\widetilde{g}^{-1}\ast
\widetilde{g}^{-1}\ast f  \ast \widetilde{A} \ast  \widetilde{A}
\ast \widetilde{\nabla}\widetilde{A}+ \widetilde{g}^{-1}\ast f \ast
\widetilde{A} \ast  \widetilde{A} \ast
 \widetilde{\nabla}\widetilde{A} \nonumber\\
 &+&S \ast \widetilde{A}\ast
\widetilde{\nabla} \widetilde{A}+A \ast S \ast
\widetilde{\nabla}\widetilde{A}+A \ast A \ast U.\nonumber
\end{eqnarray*}
 (Here $V \ast W$ denotes a linear combination of
contractions of the tensor fields $V$ and $W$ by the metric $g$.)

\hspace *{0.4cm}

Proof. As in [K1], it is easy to verify that

\begin{displaymath}
\widetilde{g}^{-1}-g^{-1}=\widetilde{g}^{-1} \ast f,
\end{displaymath}

\begin{displaymath}\nabla f=\widetilde{g} \ast P,\end{displaymath}

\begin{displaymath} \nabla
\widetilde{g}^{-1}=(\nabla-\widetilde{\nabla})\widetilde{g}^{-1}=\widetilde{g}^{-1}
\ast P,\end{displaymath}

\begin{displaymath}
 \widetilde{\nabla}W=\nabla W+P \ast W \end{displaymath}
for any tensor field $W$,

\begin{displaymath}
\widetilde{\Delta}\widetilde{A}=\Delta \widetilde{A}+f \ast
\widetilde{g}^{-1}\ast \widetilde{\nabla}^2\widetilde{A}+P \ast
\widetilde{\nabla}\widetilde{A}+Q \ast \widetilde{A}+P \ast P \ast
\widetilde{A},\end{displaymath}  and

\begin{displaymath}
\widetilde{\Delta} \widetilde{\nabla} \widetilde{A}=\Delta
 \widetilde{\nabla}  \widetilde{A}+f \ast \widetilde{g}^{-1}\ast
\widetilde{\nabla}^3 \widetilde{A}+P \ast
\widetilde{\nabla}^2\widetilde{A}+Q \ast
\widetilde{\nabla}\widetilde{A}+P*P*\widetilde{\nabla}\widetilde{A}. \end{displaymath}

Recall also that

\begin{displaymath}
\frac{\partial}{\partial t}\nabla P=\nabla \frac{\partial P}{\partial t}+\frac{\partial \Gamma}{\partial t}*P.
\end{displaymath}
Then Lemma 2 follows from Lemma 1 by direct computations.         \hfill{$\Box$}

\hspace *{0.4cm}

Now  as in [K1] we let
\begin{displaymath}
\mathcal{X}=T_2(M) \bigoplus T_3(M), \hspace{2mm}  \mathcal{Y}=T_2(M)\bigoplus
T^1_2(M)\bigoplus T^1_3(M),
\end{displaymath}
\noindent and for each $t\in [0,T]$ let
\begin{displaymath}
\mathbf{X}(t)=S(t)\bigoplus U(t) \in \mathcal{X},  \hspace{2mm}
\mathbf{Y}(t)=f(t)\bigoplus P(t)\bigoplus Q(t) \in \mathcal{Y},
\end{displaymath}
\noindent where $S, U, f, P$ and $Q$ are defined as above.

\noindent Then we have the following

\hspace *{0.4cm}

{\bf Lemma 3} Assume that the manifold $M^n$ is orientable. Let $F_t, \widetilde{F}_t:M^n \rightarrow \mathbb{R}^{n+1}$ be
two complete solutions of the mean curvature flow on $M^n \times
[0,T]$ with $|A|_{g_t}\leq K$ and $|\widetilde{A}|_{\widetilde{g}_t}
\leq \widetilde{K}$ for some constants $K$ and $\widetilde{K}$.
Suppose $F_T=\widetilde{F}_T$. Then for any $0< \delta <T$, there
exists a positive constant $C=C(\delta, K,\widetilde{K},T)$ such
that
\begin{eqnarray*}
|(\frac{\partial}{\partial t}-\Delta_{g_t})\mathbf{X}|^2_{g_t}
&\leq & C(|\mathbf{X}|^2_{g_t}+|\mathbf{Y}|^2_{g_t}),\\
|\frac{\partial}{\partial t}\mathbf{Y}|^2_{g_t} &\leq &
C(|\mathbf{X}|^2_{g_t}+ |\nabla
\mathbf{X}|^2_{g_t}+|\mathbf{Y}|^2_{g_t}).
\end{eqnarray*}

\hspace *{0.4cm}

Proof. By Ecker-Huisken [EH] there exist constants $C_m=C_m(\delta,
K,T)$ and $\widetilde{C}_m=\widetilde{C}_m(\delta, \widetilde{K},T)$
such that $|\nabla^mA|_{g_t}\leq C_m$ and $|\widetilde{\nabla}^m
\widetilde{A}|_{\widetilde{g}_t}\leq \widetilde{C}_m$ on $M^n \times
[\delta,T]$.

 Since $|A|_{g_t}\leq K$, it follows from Lemma 1 (2.1) that the  metrics $\{g_t \}_{t \in [0,T]}$ are uniformly equivalent.
 Similarly, the metrics $\{ \widetilde{g}_t\}_{t \in [0,T]}$
 are uniformly equivalent too.  But  by our assumption $F_T=\widetilde{F}_T$, and $g_T=\widetilde{g}_T$, so
 $\{g_t \}_{t \in [0,T]}$ and $\{ \widetilde{g}_t\}_{t \in [0,T]}$
 are equivalent to each other. It follows that $|\widetilde{g}^{-1}|_{g_t}$,$|\widetilde{\nabla}^m \widetilde{A}|_{g_t}$,
  $|f|_{g_t}$, $|S|_{g_t}$, and $|U|_{g_t}$ are bounded.

  Now we see  that $|P|_{g_t}$ is bounded by using the second formula in Lemma 2 and the assumption $P(T)=0$. In
  fact, for any $x \in M^n$,
\begin{displaymath}
|P(x,t)|_{g_t}=|P(x,T)-P(x,t)|_{g_t} \leq  \int^T_t|\frac{\partial
P}{\partial t}(x,s)|_{g_t}ds \leq C'.
\end{displaymath}
(One can also prove this using Lemma 1 (2.2). Compare with [K1].)

   Similarly $Q$ and $\nabla^m P$ are bounded.
  Then Lemma 3 follows from Lemma 2.      \hfill{$\Box$}

\hspace *{0.4cm}

Now as above, let  $F_t,
\widetilde{F}_t:M^n \rightarrow \mathbb{R}^{n+1}$ be two complete
solutions of the mean curvature flow on $M^n \times [0,T]$ with
bounded second fundamental forms, where $M^n$ is connected and orientable. Suppose $F_T=\widetilde{F}_T$.

Using the identity
\begin{displaymath}\nabla^m\widetilde{\nabla}^l\widetilde{A}=\nabla^{m-1}\widetilde{\nabla}^{l+1}\widetilde{A}+\sum\limits_{i=0}^{m-1}\nabla^iP
\ast
\nabla^{m-1-i}\widetilde{\nabla}^l\widetilde{A}\end{displaymath} one
sees that $\nabla S=\nabla A-\nabla \widetilde{A}$ and $\nabla
U=\nabla^2 A-\nabla \widetilde{\nabla}  \widetilde{A}$ are bounded
on $M^n \times [\delta, T] $ for any $0< \delta <T$. So the required
growth condition of [K1,Theorem 3.1] is  verified.

With the help of Lemma 3,  we can apply  [K1,Theorem 3.1] to conclude
that $\mathbf{X}=0$, $\mathbf{Y}=0$ on $M^n \times [\delta, T]$ for
any $0< \delta <T$. Then by the uniqueness theorem for hypersurfaces in a Euclidean space (see for example Theorem 6.4 in Chapter VII of [KN]), for each $t\in [\delta,T]$, $F_t$ and $\widetilde{F}_t$ coincide up to an ambient Euclidean isometry. In particular, there exists a Euclidean isometry $\bar{\sigma}: \mathbb{R}^{n+1}\rightarrow \mathbb{R}^{n+1} $ such that $\bar{\sigma}\circ F_\delta=\widetilde{F}_\delta$.

Now $\bar{\sigma}\circ F_t$ and $\widetilde{F}_t$ are two complete solutions of the mean curvature flow on $M^n \times [\delta, T]$ with
bounded second fundamental forms and with the same initial value. By Chen-Yin's uniqueness theorem for the mean curvature flow [CY], $\bar{\sigma}\circ F_t=\widetilde{F}_t$ for any $t\in [\delta,T]$. In particular,
$\bar{\sigma}\circ F_T=\widetilde{F}_T$. Combining  with our assumption we get $\bar{\sigma}\circ F_T= F_T$. It follows that either $\bar{\sigma}=Id$ or the image of $F_T$ is a hyperplane in $\mathbb{R}^{n+1}$ and $\bar{\sigma}$ is a reflection w.r.t. it.  In the latter case, by  using what we have proved in the previous  paragraph with $\widetilde{F}_t$ there replaced by the trivial hyperplane solution to the mean curvature flow, we see that the image of $F_t$ is also a hyperplane for any $t\in [\delta,T]$.   So in both cases $F_t=\widetilde{F}_t$ for any $t\in [\delta,T]$. Since $\delta \in (0,T)$ can be arbitrarily small, by continuity the Theorem is proved.

\hspace *{0.4cm}

{\bf Acknowledgements}  I was partially supported by NSFC no.10671018 and by  Laboratory of Mathematics and Complex Systems, Ministry of Education, at BNU. The first version of this note was  posted on the arXiv in 2009. Recently there appeared two papers extending our result above to the higher codimension case, cf. [LM] (where the ambient spaces may be certain general Riemannian manifolds) and [Z] (where the ambient spaces are Euclidean).  I would like to thank the authors of these two papers for their comments on the first version of my note, in particular, thank Man-Chun Lee for pointing out a gap in the  argument in it.
I would also like to thank the referee for the comments.

\bibliographystyle{amsplain}

\hspace *{0.4cm}

{\bf References}

[CY] B.-L. Chen, L. Yin, Uniqueness and pseudolocality theorems of the mean curvature
flow, Comm. Anal. Geom. 15 (2007), no. 3, 435-490.

[EH] K. Ecker, G. Huisken, Interior estimates for hypersurfaces moving by mean curvature,
Invent. Math. 105 (1991), 547-569.

[H] G. Huisken, Flow by mean curvature of convex
surfaces into spheres, J. Diff. Geom. 20 (1984), no. 1, 237-266.

[KN] S. Kobayashi,  K. Nomizu,  Foundations of differential geometry. Vol. II.  John Wiley \& Sons, Inc., 1996.

[K1] B. Kotschwar, Backwards uniqueness of the Ricci flow, Inter.  Math. Res. Not. 2010, no. 21, 4064-4097.

[K2] B. Kotschwar, A short proof of backward uniqueness for some geometric evolution equations,
Internat. J. Math. 27 (2016), no. 12, 1650102, 17 pp.

[LM] M.-C. Lee, M. S. John Ma,
Uniqueness Theorem for non-compact mean curvature flow with possibly unbounded curvatures,  arXiv:1709.00253.

[Z] Z.-H. Zhang,
A note on the backwards uniqueness of mean curvature flow, Sci. China Math. (2018). https://doi.org/10.1007/s11425-017-9231-4. arXiv:1709.00798.

\end{document}